\documentclass[12pt]{amsart}
\usepackage[utf8]{inputenc}

\usepackage{amsrefs}
\usepackage{supertabular}
\usepackage{amsxtra,amssymb,amsmath,amscd,url,enumitem}

\setenumerate{itemsep=5pt}
\setitemize{itemsep=5pt}

\newlist{enumth}{enumerate}{1}
\setlist[enumth]{label=\emph{(\arabic*)}, ref=(\arabic*)}

\usepackage[utf8]{inputenc}
\usepackage[english]{babel}
\usepackage{fullpage}
\usepackage{scrtime}
\usepackage[colorlinks]{hyperref}
\usepackage{datetime}
\usepackage[all]{xy}
\usepackage{tikz-cd}
\usepackage{mathrsfs}
\usepackage{xspace}

\usepackage{tensor}
\usepackage{braket}

\usepackage{todonotes}
\usepackage{tensor}

\shortdate

\DeclareMathSymbol{A}{\mathalpha}{operators}{`A}%
\DeclareMathSymbol{B}{\mathalpha}{operators}{`B}%
\DeclareMathSymbol{C}{\mathalpha}{operators}{`C}%
\DeclareMathSymbol{D}{\mathalpha}{operators}{`D}%
\DeclareMathSymbol{E}{\mathalpha}{operators}{`E}%
\DeclareMathSymbol{F}{\mathalpha}{operators}{`F}%
\DeclareMathSymbol{G}{\mathalpha}{operators}{`G}%
\DeclareMathSymbol{H}{\mathalpha}{operators}{`H}%
\DeclareMathSymbol{I}{\mathalpha}{operators}{`I}%
\DeclareMathSymbol{J}{\mathalpha}{operators}{`J}%
\DeclareMathSymbol{K}{\mathalpha}{operators}{`K}%
\DeclareMathSymbol{L}{\mathalpha}{operators}{`L}%
\DeclareMathSymbol{M}{\mathalpha}{operators}{`M}%
\DeclareMathSymbol{N}{\mathalpha}{operators}{`N}%
\DeclareMathSymbol{O}{\mathalpha}{operators}{`O}%
\DeclareMathSymbol{P}{\mathalpha}{operators}{`P}%
\DeclareMathSymbol{Q}{\mathalpha}{operators}{`Q}%
\DeclareMathSymbol{R}{\mathalpha}{operators}{`R}%
\DeclareMathSymbol{S}{\mathalpha}{operators}{`S}%
\DeclareMathSymbol{T}{\mathalpha}{operators}{`T}%
\DeclareMathSymbol{U}{\mathalpha}{operators}{`U}%
\DeclareMathSymbol{V}{\mathalpha}{operators}{`V}%
\DeclareMathSymbol{W}{\mathalpha}{operators}{`W}%
\DeclareMathSymbol{X}{\mathalpha}{operators}{`X}%
\DeclareMathSymbol{Y}{\mathalpha}{operators}{`Y}%
\DeclareMathSymbol{Z}{\mathalpha}{operators}{`Z}%

\renewcommand{\leq}{\leqslant}
\renewcommand{\geq}{\geqslant}
 
\numberwithin{equation}{section}

\newcommand{\aj}[1]{(N-\lfloor {{#1}}/2\rfloor)(N-\lceil{{#1}}/2\rceil)}
\newcommand{\bj}[1]{N(N-{{#1}})}

\renewcommand{\mathcal}{\mathscr}

\newcommand{\Zz}{\mathbf{Z}}

\newcommand{\Rr}{\mathbf{R}}

\def\loccit{loc.\kern3pt cit.{}\xspace}
\def\cf{see\kern.3em}
\def\Cf{See\kern.3em}
\def\eg{e.g.\kern.3em}

\def\resp{\text{resp.}\kern.3em}

\renewcommand{\rho}{\varrho}

\DeclareMathSymbol{\gena}{\mathord}{letters}{"3C}
\DeclareMathSymbol{\genb}{\mathord}{letters}{"3E}

\theoremstyle{plain}
\newtheorem{theorem}{Theorem}[section]
\newtheorem*{theorem*}{Theorem}
\newtheorem{lemma}[theorem]{Lemma}
\newtheorem{corollary}[theorem]{Corollary}

\newtheorem{proposition}[theorem]{Proposition}

\theoremstyle{remark}

\theoremstyle{definition}

\newtheorem{definition}[theorem]{Definition}

\newtheorem{example}[theorem]{Example}
\newtheorem{remark}[theorem]{Remark}

\renewcommand{\geq}{\geqslant}
\renewcommand{\leq}{\leqslant}

\newcommand{\leqa}{\leq_{\mathrm{as}}}

\newcommand{\lea}{<_{\mathrm{as}}}
\newcommand{\gea}{>_{\mathrm{as}}}

\begin{document}

\title{The gaps in the multiplication table}

\author{Emmanuel Kowalski}
\address[E. Kowalski]{D-MATH, ETH Z\"urich, R\"amistrasse 101, 8092 Z\"urich, Switzerland} 
\email{kowalski@math.ethz.ch}

\author{Vivian Z. Kuperberg}
\address[V.Z. Kuperberg]{D-MATH, ETH Z\"urich, R\"amistrasse 101, 8092
  Z\"urich, Switzerland}
\email{vivian.kuperberg@math.ethz.ch}

\begin{abstract}
  We find the explicit list of the gaps between successive elements of
  the $N\times N$ multiplication table for any positive integer~$N$. 
\end{abstract}

\makeatletter
\@namedef{subjclassname@2020}{%
  \textup{2020} Mathematics Subject Classification}
\makeatother



\maketitle

\begin{flushright}
  \textit{Dedicated to Zeev Rudnick, an amateur of fine gaps}
\end{flushright}

\section{Introduction}


This paper is about the following simple-sounding question: \emph{What
  are the gaps which occur between successive elements of the
  multiplication table listing all the products $ab$ for $a$ and $b$
  integers between~$1$ and~$N$?} Readers may be interested in this
problem because they like the multiplication table, or because they
are fascinated by the general problem of understanding the gaps
between elements of various sequences.  Before stating the result, we
survey some of the mysteries of both of these areas of mathematics.

\subsection{The multiplication table}

The multiplication table is a humble, yet very famous, finite set of
integers: in fact probably the best known set for a large fraction of
the world's population! It turns out, as Erd\H os observed in the
1930's, that this set has remarkable properties. The key observation
is that the set of values $ab$ with $1\leq a,b\leq N$ is unusual for
the simple reason that its elements typically have \emph{twice as
  many} prime factors as most integers between~$1$ and~$N^2$ (this is
true whether the prime factors are counted with or without
multiplicity). More precisely, if~$N$ is large then, according to
results of Hardy and Ramanujan (later superseded by the Erd\H os--Kac
Theorem, a simple proof of which is explained for instance
in~\cite[\S\,2.3]{pnt}), ``most'' integers~$n\leq N$ have about
$\log(\log(N))$ prime factors (either with or without
multiplicity). Since~$\log(\log(N^2))$ differs little from
$\log(\log(N))$, the same holds for most integers~$n\leq N^2$. By
contrast, a typical number $ab$ in the multiplication table has about
$2\log(\log(N))$ prime factors.  In particular, it follows fairly
easily that the number of distinct elements in the multiplication
table, say $M_N$, is $o(N^2)$, i.e., the ratio $M_N/N^2$ tends to~$0$
as $N\to +\infty$ (see also the paper of Koukoulopoulos and
Thiel~\cite{k-t}).

Can one be more precise? This is known as the \emph{multiplication table
  problem}, and many mathematicians investigated the size of~$M_N$
(among whom we cite in particular Tenenbaum~\cite{tenenbaum}) before
Ford~\cite[Cor.\,3]{ford} obtained the precise order of magnitude. His
striking statement (which had been conjectured by Erd\H os in 1960) is
that if we set
\[
  \alpha=\int_1^{1/\log(2)}\log(t)dt=
  1-\frac{1+\log(\log(2))}{\log(2)}=0.08607\ldots,
\]
then there exist constants $c_1$, $c_2>0$ such that, for all
integers~$N\geq 3$, the inequalities
\begin{equation}\label{eq-ford}
  c_1(\log N)^{-\alpha}(\log\log N)^{-3/2} \leq \frac{M_N}{N^2}\leq c_2(\log
  N)^{-\alpha}(\log\log N)^{-3/2}
\end{equation}
hold.

The multiplication table problem appears, surprisingly, in other
applications. One of them is found in the paper~\cite{k-t}, previously
mentioned. Another is due to Brent and Kung~\cite[\S\,3.2]{b-k}, who
rediscovered the problem while studying a very practical question in
computer science. More recently, Eberhard, Ford and Green~\cite{efg}
proved an analog of Ford's theorem for \emph{random permutations},
where the same constant~$\alpha$ appears (see also Granville's
entertaining survey of such analogies~\cite{granville2}).

It remains unclear whether the quantity
\[
\frac{M_N}{N^2}\times (\log  N)^{\alpha}(\log\log N)^{3/2}
\]
actually converges or not; this suggests the interest of doing
numerical experiments, which has been done most recently by Brent,
Pomerance, Purdum and Webster~\cite{bppw}. (We note that Green
and Sawhney have announced progress towards understanding the limiting
behavior of this quantity in a recent Oberwolfach meeting.)

Our own work also started with numerical experiments with the full
multiplication table, although on a fairly small scale. This brings us
to the general context of the question we study.

\subsection{Gaps in sequences}

Sequences of integers or real numbers are among the most basic objects
of mathematics. In many cases, their study involved experimental
aspects, and some of the most famous results of arithmetic have been
discovered in this manner. Notable examples include lists of primes,
of course\footnote{\ This led Legendre and Gauss to guess a correct
  asymptotic formula concerning the number of primes $p$ up to a large
  number $x$; see~\cite{tschinkel} for a nice look at Gauss's
  tables.}, but also tables of sums of two squares\footnote{\ In some
  sense, much of algebraic number theory can be motivated from these
  tables. They certainly suggest the theorem of Fermat that primes
  which are $1$ modulo~$4$ are sums of two squares of
  integers--famously provable in one line~\cite{zagier}.}, or lists of
zeros of the Riemann zeta function.\footnote{\ The first few zeros
  were computed by Riemann--see for instance the discussion and
  reproduction of Riemann's papers in Edwards's
  book~\cite[p.\,156]{edwards}.}

Since the integers (and the real numbers) are ordered, it happens
without much thinking that one studies a sequence of non-negative
numbers $(a_n)_{n\geq 1}$ by ordering it, or rather by ordering its
initial terms
\[
  a_1,\ a_2,\ a_3,\ \ldots, a_N,
\]
where $N\geq 1$ will vary among positive integers. There is
some permutation of these terms which orders them, say
\[
  a_{\sigma(1)}\leq a_{\sigma(2)}\leq \cdots\leq a_{\sigma(N)}.
\]

As soon as this is done, mathematicians will almost reflexively start
looking for patterns, and among the most immediate question they will
try to understand is that of the shape of the \emph{gaps} between
successive terms: do the differences
\[
  a_{\sigma(2)}-a_{\sigma(1)},\quad\ldots,\quad
  a_{\sigma(N)}-a_{\sigma(N-1)}
\]
exhibit any interesting features?

Even for sequences $(a_n)$ which are extremely well-understood, hidden mysteries are revealed by this new
question--many of which lead to
strikingly challenging conjectures.  Although the study of the
multiplication table is quite different, we provide some examples here.

\begin{enumerate}
\item The ``three gaps theorem'', conjectured by Steinhaus and proved
  by Sós~\cite{sos}, states that there are \emph{at most three}
  distinct gaps between the fractional parts of the numbers $\xi$,
  $2\xi$, \ldots, $N\xi$, if $\xi$ is an irrational number. The values
  of these gaps, of course, depend on~$N$ and~$\xi$, but there are
  always at most three distinct gaps! There is a recent appealing
  proof of this remarkable fact due to Marklof and Strömbergsson
  (see~\cite{m-s}).

\item Gaps between primes are one of the crowning glories of analytic
  number theory. Of particular note are the proof by
  Y. Zhang~\cite{zhang} that \emph{some} integer occurs infinitely
  often as the difference between two consecutive primes, and the
  results of J. Maynard~\cite{maynard} which both simplified and
  strengthened Zhang's result by obtaining bounded intervals containing
  any fixed number of primes.

\item The gaps between ordinates of the zeros of the Riemann zeta
  function were famously studied by
  Montgomery~\cite{montgomery}. Thanks to a chance encounter with
  Dyson, this led to the discovery of links between the zeros of the
  zeta function (and later its generalizations) and eigenvalues of
  random unitary matrices, as recounted for instance in the survey of
  Katz and Sarnak~\cite{katz-sarnak}.

\item The energy levels of quantum mechanical systems frequently form
  an infinite sequence tending to $\infty$, and the gaps between these
  energy levels have considerable importance in physics and
  mathematical physics. The eigenvalues of the Laplace operator on
  compact Riemannian manifolds provide examples of such
  sequences. Maybe the simplest concrete cases are the sequences of
  values of a positive definite quadratic form corresponding to the
  Laplace operator on a torus (such as $\alpha n^2+\beta m^2$ for
  positive numbers $\alpha$ and~$\beta$). Understanding the
  distribution of these gaps is already an extraordinarily challenging
  problem with connections to many areas of mathematics (see Sarnak's
  survey~\cite{sarnak} for a general discussion).
\end{enumerate}

\subsection{Gaps in the multiplication table}

The goal of this paper is to examine the gaps between successive
elements in the multiplication table of the first $N$ integers,
allowing~$N$ to be arbitrarily large, thus combining the two topics
discussed previously.

This is a good example, because although computing the integers in the
$N\times N$ multiplication table, for small enough~$N$, is completely
straightforward, it is also true that rearranging them in order will
scramble these numbers quite a bit, with un unpredictable result.
Taking $N=5$, we go from
\[
  1, 2, 3, 4, 5, 2, 4, 6, 8, 10, 3, 6, 9, 12, 15, 4, 8, 12, 16, 20, 5,
  10, 15, 20, 25
\]
to
\[
  1, 2, 2, 3, 3, 4, 4, 4, 5, 5, 6, 6, 8, 8, 9, 10, 10, 12, 12, 15, 15,
  16, 20, 20, 25,
\]
with successive gaps
\[
  1, 0, 1, 0, 1, 0, 0, 1, 0, 1, 0, 2, 0, 1, 1, 0, 2, 0, 3, 0, 1, 4, 0,
  5.
\]

As in this example, we will keep track of the multiplicities, both in
the multiplication table and in the gaps (there may be many pairs of
successive elements of the multiplication table differing by~$1$, for
instance). In particular, we will be interested in gaps which occur
exactly once, which we call ``isolated gaps''.  We will sometimes call
``gap'' any difference $cd-ab$ of elements $ab\leq cd$ of the
multiplication table, and call ``primitive gap'' those differences
which really occur between consecutive elements. In other words, primitive gaps
are the gaps between $ab$ and $cd$ for which the inequalities
$$
ab\leq ef\leq cd
$$
with $1\leq e,f\leq N$ imply that $ef=ab$ or $ef=cd$. Moreover, we
will sometimes denote the gap $cd-ab$ (primitive or not) by stating an
inequality
\[
  ab<cd
\]
with $1\leq a,b,c,d\leq N$, and say that the gap ``occurs from the
interval'' $ab<cd$.

Our main result \emph{exactly} determines the primitive gaps in the
multiplication table. To motivate the result, and suggest how one can
discover it, we start with some examples which demonstrate the
surprising patterns that occur.

We take $N=100$. The following is a list $(m_0,m_1, \ldots)$ of
non-negative integers where~$m_{i}$ is the number of primitive gaps
equal to~$i$ when listing in increasing order, with multiplicity, the
elements of the $100\times 100$ multiplication table, with the list
ending at the largest gap:
\begin{gather*}
  (7094, 1039, 656, 380, 265, 142, 124, 64, 60, 36, 40, 22, 16, 10,
  11, 10, 6, 5, 1, 1, \\
  0, 0, 0,  0, 0, 0, 0, 0, 1, 0, 0, 0, 0, 0, 0, 0, 1, 0, 0, 0, 0, 0,
  0, 0, 1, 0, 0, 0, 0, 0, 0, 1, 0, 0,
  \\
  0, 0, 0, 0, 1, 0, 0, 0, 0, 0, 1, 
  0, 0, 0, 0, 0, 1, 0, 0, 0, 0, 1, 0, 0, 0, 0, 1, 0,
  0, 0, 1, 0, 0,\\
  0, 1, 0, 0, 1, 0, 0, 1, 0, 1, 0, 1, 1, 1).
\end{gather*}
For instance, there are $1039$ pairs $(m,m+1)$
which belong to the multiplication table.

We see some patterns emerging: first, there is an interval of initial
integers that are all primitive gaps, which is then followed by a
sequence of more sparse \emph{isolated} primitive gaps. Moreover, if we
count the gaps, we see that the largest gap is $100$ exactly.

We can very quickly check that this last fact is general: among the
elements of the $N\times N$ multiplication table, the largest gap
between successive elements is $N$. Indeed, if $ab$ is an element of the
table, with $1\leq a\leq b\leq N$, and if furthermore $a<N$, then
$(a+1)b$ is also in the $N\times N$ multiplication table, and we have
\[
  (a+1)b-ab=b\leq N.
\]

These basic features of the list of gaps turn out to be general facts.
But looking at other examples suggests rather more. For instance, here
is the list of gaps in the $105\times 105$ multiplication table,
presented in the same way:
\begin{gather*}
  (7795, 1171, 740, 421, 292, 141, 138, 66, 71, 32, 36, 28, 21, 11, 13,
  12, 9, 5, 4, \\
  0, 0, 0, 0, 0, 1, 0, 0, 0, 0, 0, 0, 0, 0, 1, 0, 0, 0, 0, 0, 0, 0, 1,
  0, 0, 0, 0, 0, 0, 0, 1, 0, 0, 0, 
  \\
  0, 0, 0, 1, 0, 0, 0, 0, 0, 0, 1, 0, 0, 0, 0, 0, 1, 0, 0, 0, 0, 0, 1, 0,
  0, 0, 0, 1, 0, 0, 0, 0, 1, 0, 
  \\
  0, 0, 1, 0, 0, 0, 1, 0, 0, 1, 0, 0, 1, 0, 1, 0, 1, 1, 1).
\end{gather*}

We see the pattern of an interval followed by isolated gaps.  But more
importantly, if we \emph{reverse} the list by starting from the end
(which we know corresponds to a gap equal to $N$), we find
\begin{gather*}
  N=100:\ 1, 1, 1, 0, 1, 0, 1, 0, 0, 1, 0, 0, 1, 0, 0, 0, 1, 0, 0, 0, 1,
  0, 0, 0, 0, 1, 0, 0, \ldots 
  \\
  N=105:\ 1, 1, 1, 0, 1, 0, 1, 0, 0, 1, 0, 0, 1, 0, 0, 0, 1, 0, 0, 0, 1,
  0, 0, 0, 0, 1, 0, 0, \ldots 
\end{gather*}

After reversion, the ``largest'' gaps in the multiplication
table are (apparently) independent of $N$.  Next, we list the
integers~$k\geq 0$ indicating where the $1$'s in this common list
appear (which means that the corresponding gaps are $N-k$); we get
\[
  k=0,1,2,4,6,9,12,16,\ldots
\]

We will prove that this patterns holds in general, and identify the
sequence of values~$k$ above.
Here is the precise result.

\begin{theorem}\label{th-table}
  Let $N\geq 2$ be an integer and let $R$ be the largest integer~$j$
  such that~$\lfloor j^2/4\rfloor<N$.  Define the integer
  \[
    T=\begin{cases}
      R&\text{ if $N$ is of the form $k^2$ or $k^2+k$,}
      \\
      R-1&\text{ if $N$ is of the form $k^2-1$ or~$k^2+k-1$,}
      \\
      R-2&\text{ otherwise},
    \end{cases}
  \]
  where in the first two cases, $k$ is a positive integer.

  Then the set of primitive gaps in the $N\times N$ multiplication table
  is the union of the sets $\mathcal{S}_N$ and $\mathcal{B}_N$ defined
  by
  \[
    \mathcal{S}_N=\{0,\ldots, T\}\text{ and } \mathcal{B}_N=\{N-\lfloor
    j^2/4\rfloor
    \,\mid\, 1\leq j\leq R-2\}.
  \]
  Moreover, the two sets are disjoint, and in fact for
  all~$x\in\mathcal{S}_N$ and~$y\in\mathcal{B}_N$, we have $x<y$.
\end{theorem}

\begin{remark}
  (1) The integer~$R$ is of size roughly $2\sqrt{N}$ (see
  Lemma~\ref{lm-r} below for a more precise description). Thus we have
  about $4\sqrt{N}$ distinct primitive gaps among members of the
  multiplication table.
  
  (2) Theorem~\ref{th-table} is a new illustration of how far the
  multiplication table is from a ``generic'' set of integers. Indeed,
  for a random set $A\subset \{1,\ldots,N^2\}$ with relatively large
  density $\delta=|A|/N^2$, one would expect that every interval of
  length somewhat larger than this should contain some element of~$A$
  (this is the same basic reasoning behind the predictions of the
  ``Cramér model'' for the size of gaps between primes; see for
  instance Granville's elementary discussion~\cite{granville} of this
  model). In particular, it should also be the case that there are at
  most about $1/\delta$ distinct primitive gaps between elements
  of~$A$.  In the case of the multiplication table, Ford's
  result~(\ref{eq-ford}) gives $\delta$ of size about
  $(\log N)^{-0.8}$, so it is striking that there are actually about
  $4\sqrt{N}$ distinct gaps.

  Thus Theorem~\ref{th-table} shows that there are much larger gaps than
  expected on probabilistic grounds. 
  It is usually extremely challenging to prove that a given number is
  a primitive gap in a sequence, and it is even more so when there are
  so many possible gaps. The point, as will be seen in the proof of
  Theorem~\ref{th-table}, is that the gaps, even the small ones, are
  located among very large values of the multplication table, and
  arise from some type of ``algebraic'' rigidity.

  (3) One would expect that the bulk of the multiplication table is
  more random.  One can show that an interval of the form $X<X+Y$,
  with $N\leq X$, $X\leq cN^{2}$ and $Y\geq c X^{1/3}$ (for suitable
  constant $c>0$) contains approximately as many elements of the
  multiplication table as expected heuristically (counted with
  multiplicity). This results from a fairly elementary argument with
  exponential sums and exponent pairs
  (see, e.g., the treatment by Graham and
  Kolesnik~\cite[\S\,4.3]{gk}).

  (4) We also note that although the theorem shows that the large gaps
  follow rigid ``algebraic'' rules, the multiplicities of the other
  gaps are certainly much more arithmetic in nature. Let $\mu_N(n)$
  denote the multiplicity of~$n\geq 0$ as a gap in the $N\times N$
  multiplication table. Let further
  $$
  d_N(n)=\sum_{\substack{ab=n\\a,b\leq N}}1
  $$
  denote the ``representation function'' for the multiplication
  table. Then 
  \[
    \mu_N(0)=
    N^2-M_N,
  \]
  where~$M_N$ is the number of distinct elements in the multiplication
  table, whose subtle behavior was already mentioned (see~(\ref{eq-ford})).

  Similarly, since any gap equal to~$1$ is primitive, we have
  \[
    \mu_N(1)
    =\sum_{\substack{1\leq a,b,c,d\leq N\\ab-cd=1}}\frac{1}{d_N(ab)d_N(cd)}.
  \]

  One knows (see the result~\cite[Th.\,1.1,\,(3.1)]{ganguly-guria} of
  Ganguly and Guria) that
  $$
  \sum_{\substack{1\leq a,b,c,d\leq N\\ab-cd=1}}1\sim
  \frac{2}{\zeta(2)}N^2,
  \quad\quad
  \zeta(2)=\sum_{n\geq 1}\frac{1}{n^2},
  $$
  as $N\to +\infty$, which implies that~$\mu_N(1)$ is roughly of size
  $N^2$. 
  It would be interesting to determine the asymptotic behavior of
  $\mu_N(k)$ for $k$ fixed, or when $k\leq N^{1/2-\delta}$ for
  some~$\delta>0$.

  (5) Experimentally, it seems that the isolated primitive gaps are
  those in the set $\mathcal{B}_N$ and the gaps~$R-1$ and~$R$ which
  occur for special values of~$k$.
  In other words, all gaps $\leq
  R-2$ seem to occur with multiplicity $\geq
  2$. For those which are $\leq
  R-5$, this fact can be proved easily from the proof of
  Proposition~\ref{pr-existence-gaps}, (2) below, but the remaining
  cases are more elusive.
\end{remark}

\begin{example}\label{ex-examples}
  We give one additional example to illustrate the main result. We list
  the multiplicity distribution of the gaps as a tuple
  $$
  (\mu_N(0),\mu_N(1),\ldots, \mu_N(N)),
  $$
  meaning that~$0$ is a primitive gap with multiplicity $\mu_N(0)$, etc.

  Let $N=73$, which is of the form $k^2+k+1$ with~$k=8$. We get
  \begin{gather*}
    (3678, 638, 377, 227, 136, 71, 63, 25, 22, 23, 16, 6, 10, 9, 10,
    2, 0, 1, 0, 0, 0, 0, 0,\\
    0, 1, 0, 0, 0, 0, 0, 0, 1, 0, 0, 0, 0, 0,
    1, 0, 0, 0, 0, 0, 1, 0, 0, 0, 0, 1, 0, 0,\\ 0, 0, 1, 0, 0, 0, 1, 0,
    0, 0, 1, 0, 0, 1, 0, 0, 1, 0, 1, 0, 1, 1, 1).
  \end{gather*}

  We have $R=17=2k+1$ and 
  $$
  \mathcal{B}_{73}=\{17,24,41,37,43,48,53,57,61,63,67,69,71,72,73\}.
  $$
  
  In this case, $R=17$ is an element (the smallest)
  of~$\mathcal{B}_{73}$; in fact, $R$ is an element of~$\mathcal B_{N}$ whenever $N=k^2+k+1$ or $N=k^2+1$
  for some integer $k\geq 1$ (this can be checked, for instance, using
  Lemma~\ref{lm-r} below).
\end{example}

\section{The ``asymptotic'' multiplication table}\label{sec-large}

The basic idea behind the proof of Theorem~\ref{th-table} is simply
that ``large'' elements of the multiplication table are ordered ``like
polynomials'', in the sense that if $a$ and $b$ are suitably small
compared to~$N$, then the values $(N-a)(N-b)$ are ordered like the
values of the polynomials $(X-a)(X-b)$ as $X$ tends to infinity. This
ordering of polynomials is very explicit and the corresponding gaps
can be determined easily.  Although we could dispense with this step,
we first look at polynomials, since it motivates the shape of the
set~$\mathcal{B}_N$ in Theorem~\ref{th-table} and was the starting
point of the investigation.  We first note that the notion of gaps,
and of primitive gaps, makes sense for subsets of any ordered abelian
group.

\begin{definition}
  The \emph{asymptotic order} on $\Rr[X]$ is the total ordering~$\leqa$
  determined by $f\leqa g$ if and only if $f(x)\leq g(x)$ for all $x$
  large enough.

  The \emph{asymptotic order} on~$\Rr^d$ is the total (pre)-ordering $\leqa$
  determined by
  $$
  (a_0,\ldots, a_{d-1})\leq (b_0,\ldots,b_{d-1})
  $$
  if and only if
  $$
  \prod_{i=0}^{d-1} (X-a_i)\leqa \prod_{i=0}^{d-1} (X-b_i).
  $$
\end{definition}

The asymptotic order on polynomials is indeed a total order since the
difference of two polynomials is a polynomial, hence is either zero,
or has the same sign at its leading coefficient for all $x$ large
enough. With this ordering, $\Rr[X]$ is an ordered abelian group (in
fact, it is an ordered ring).

On~$\Rr^d$, we only have a pre-order since any tuples $(a_i)$
and~$(b_i)$ which are permutations of each other satisfy
$(a_i)\leq (b_i)$ and $(b_i)\leq (a_i)$. Up to this issue, the order
induced on~$\Zz^d$ by the asymptotic order on~$\Rr^d$ is a
well-ordering.

\begin{lemma}
  We have $(a,b)\leqa (c,d)$ if and only if either $a+b>c+d$, or
  $a+b=c+d$ and $ab\leq cd$.
\end{lemma}

\begin{proof}
  This is a direct computation from the definition.
\end{proof}

We now determine the primitive gaps between polynomials $(X-a)(X-b)$
with $a$, $b$ non-negative integers, or equivalently between pairs
$(a,b)$ for the asymptotic order.

We first explicitly note that if $j\geq 0$ is an integer, then
$\lfloor j/2\rfloor+\lceil j/2\rceil=j$ and
$\lfloor j/2\rfloor\lceil j/2\rceil=\lfloor j^2/4\rfloor$.

\begin{proposition}
  The asymptotic well-ordering of the set of pairs of non-negative
  integers $(a,b)$ with $a\leq b$ is determined by the following
  rules.

  \begin{enumth}
  \item For any integers $j$, $k$ with $0\leq k<j$, we have
    $$
    (a,j-a)\lea (b,k-b)
    $$
    for $0\leq a\leq j/2$ and $0\leq b\leq k/2$.
  \item For any integer~$j\geq 0$, we have
    $$
    (0,j)\lea (1,j-1)\lea\cdots \lea (\lfloor j/2\rfloor,\lceil
    j/2\rceil).
    $$
  \end{enumth}
\end{proposition}

\begin{proof}
  This is a straightforward consequence of the previous lemma.
\end{proof}

\begin{example}
  In decreasing order, we obtain
  $$
  (0,0)\gea(0,1)\gea(1,1)\gea(0,2)\gea(1,2)\gea(0,3)\gea(2,2)\gea(1,3),
  $$
  corresponding to the polynomial comparisons
  \begin{multline*}
    X^2\gea X(X-1)\gea (X-1)^2\gea X(X-2)\gea (X-1)(X-2)
    \\
    \gea X(X-3)\gea
    (X-2)^2\gea (X-1)(X-3),
  \end{multline*}
  valid for large values of the variable.  The successive gaps between
  these polynomials are
  $$
  X,\quad X-1,\quad 1,\quad X-2,\quad 2,\quad X-4,\quad 1.
  $$
\end{example}

\begin{corollary}
  The primitive gaps for the asymptotic order between the polynomials of
  the form $(X-a)(X-b)$ with $(a,b)$ non-negative integers such that
  $a\leq b$ are either positive integers or polynomials
  $$
  X-\lfloor j/2\rfloor \lceil j/2\rceil=X-\lfloor j^2/4\rfloor
  $$
  for some $j\geq 1$.

  More precisely, each of the above polynomial is a unique gap, and
  any positive integer $k\geq 1$ occurs exactly between pairs
  $(a,j-a)$ and $(a+1,j-a-1)$ with $j-2a-1=k$.
\end{corollary}

\begin{proof}  
  By the proposition, the gaps which occur (and which are primitive in
  view of the well-ordering) between these polynomials are either
  ``large'' gaps between successive elements of the form
  $$
  (\lfloor j/2\rfloor,\lceil j/2\rceil)\lea (j-1,0)
  $$
  for $j\geq 1$, or ``small'' gaps between successive elements of the
  form
  $$
  (a,j-a)\lea (a+1,j-a-1)
  $$
  with $j\geq 2$ and $0\leq a\leq \lfloor j/2\rfloor$. These two gaps
  are, respectively, equal to
  $$
  X(X-j+1)-(X-\lfloor j/2\rfloor)(X-\lceil j/2\rceil)=X-\lfloor
  j/2\rfloor \lceil j/2\rceil
  $$
  and
  \begin{multline*}
    (X-(a+1))(X-(j-a-1))-(X-a)(X-(j-a))\\
    =(a+1)(j-a-1)-a(j-a)=j-2a-1,
  \end{multline*}
  which implies the result.
\end{proof}

If we go back to the multiplication table of integers, we note that
for a given integer~$N\geq 1$, we will have
\[
  (N-a)(N-b)>(N-c)(N-d)
\]
if and only if
\[
  (a,b)\gea (c,d),
\]
provided $(a,b,c,d)$ are non-negative integers $(a,b,c,d)$ such that
$|ab-cd|<N$ (this is an elementary argument).  This means that the
gaps in the multiplication table will reflect those between
polynomials in certain ranges.  However, the restriction $|ab-cd|<N$
means that the actual proof of Theorem~\ref{th-table} must be more
involved.

\section{Disjointness}

We now begin the proof of Theorem~\ref{th-table}.

We first check that the elements of $\mathcal{S}_N$ are strictly
smaller than those of $\mathcal{B}_N$.

In fact, assume that $R\geq 3$, so that $\mathcal{B}_N$ is not
empty. Let~$m=N-\lfloor \frac{(R-2)^2}{4}\rfloor$ be the smallest
element of~$\mathcal{B}_N$. Then
$$
\Bigl\lfloor \frac{(R-2)^2}{4}\Bigr\rfloor
=\Bigl\lfloor \frac{R^2}{4}\Bigr\rfloor-R+1
\leq N-1-R+1=N-R,
$$
(by definition of~$R$), so that $m\geq R$. On the other hand, the
maximum of~$\mathcal{S}_N$ is~$<R$ unless~$N$ is of the form~$k^2$
or~$k^2+k$.  We compute the values of~$m$ and~$R$ in these cases,
namely
\begin{gather*}
  N=k^2,\quad R=2k-1,\quad \frac{(R-2)^2}{4}=k^2-3k+2+\frac{1}{4},\quad
  m=3k-2,\\
  N=k^2+k,\quad R=2k,\quad \frac{(R-2)^2}{4}=k^2-2k+1,\quad
  m=3k-1,
\end{gather*}
and see that the inequality~$m>R$ is always valid.

\section{Finding the gaps}\label{sec-finding}

In this section, we will prove the existence component of the first
assertion of Theorem~\ref{th-table}, i.e., prove that all the indicated
integers occur as gaps in the multiplication table.

We fix an integer~$N\geq 2$ and we denote by~$R$ the integer in
Theorem~\ref{th-table}, i.e., the largest integer~$j\geq 1$ such that
$\lfloor j^2/4\rfloor<N$. By ``multiplication table'', we will always
mean the $N\times N$ table.

Motivated by Section~\ref{sec-large}, we define the integers
\begin{align*}
a_j&=(N-\lfloor j/2\rfloor)(N-\lceil j/2\rceil)=N^2-jN+\lfloor
j^2/4\rfloor,
\\
b_j&=N(N-j)=N^2-jN
\end{align*}
for all integers~$j\geq 0$. As long as $0\leq j\leq N-1$, the integers
$a_j$ and $b_j$ belong to the $N\times N$ multiplication table.

The inequalities
\begin{equation}\label{eq-ajbj}
  b_R<a_R<b_{R-1}<\cdots < b_2<a_2<b_1=a_1<b_0=N^2
\end{equation}
are then valid, with 
(not necessarily primitive) gaps
\begin{equation}\label{eq-differences}
  a_j-b_j=\lfloor j^2/4\rfloor,\quad\quad b_{j-1}-a_j=N-\lfloor
  j^2/4\rfloor
\end{equation}
for $1\leq j\leq R$.  

We will now investigate which elements of the multiplication table
might lie in the intervals defined by the
inequalities~(\ref{eq-ajbj}). For orientation, we note that $b_R$ is
close to $N^2-2N^{3/2}$, so all results here concern the very large
elements of the multiplication table.

We will frequently use the following simple lemma.

\begin{lemma}\label{lm-parabola}
  Let~$k$ be a non-negative integer. As $(a,b)$ vary over pairs of
  non-negative integers with $a\leq b$ and~$a+b=k$, we have
  \begin{equation}\label{eq-ab}
    ab\leq \Bigl\lfloor\frac{k^2}{4}\Bigr\rfloor,
  \end{equation}
  with equality if and only if $a=\lfloor k/2\rfloor$ and
  $b=\lceil k/2\rceil$ and moreover either~$a=0$ or
  \begin{equation}\label{eq-ab1}
    ab\geq k-1,
  \end{equation}
  with equality if and only if~$a=1$.
\end{lemma}

\begin{proof}
  This is elementary calculus.
\end{proof}

We begin with the interval between $b_j$ and $a_j$. Given arbitrary
integers $(a,b)$ and $j\geq 0$ with $\lfloor j^2/4\rfloor\leq N$, we
have
\begin{equation}\label{ineq-1}
b_j\leq (N-a)(N-b)\leq a_j
\end{equation}
if and only if
\begin{equation}\label{ineq-1b}
0\leq (j-(a+b))N+ab\leq \lfloor j^2/4\rfloor.
\end{equation}

\begin{lemma}\label{lm-short}
  Assume that $N\geq 8$, that $0\leq a\leq b\leq N$ and that
  $\lfloor j^2/4\rfloor<N$.
  \par
  \begin{enumth}
  \item If $a+b=j$ then~\emph{(\ref{ineq-1})} holds.
    \par
  \item If~\emph{(\ref{ineq-1})} holds and $a+b\not=j$, then $a+b=j+1$
    and the inequalities
    \begin{gather*}
      N\leq ab\leq N+\frac{3}{2}\sqrt{N},\\
      b_j=N(N-j)\leq (N-a)(N-b)\leq b_j+\frac{3}{2}\sqrt{N}
    \end{gather*}
    hold.
  \end{enumth}
\end{lemma}

\begin{proof}
  The assumption on $j$ implies that $j^2/4\leq N+1$, so
  $j\leq 2\sqrt{N+1}$.
  
  If $a+b=j$, then it is elementary that~(\ref{ineq-1b}) holds, hence
  also~(\ref{ineq-1}). We assume therefore that this is not the case.

  Suppose first that $ab<N$. Then the inequality $b_j\leq (N-a)(N-b)$
  holds if and only if $a+b=j$ or $j>a+b$. In the latter case, we have
  $$
  (j-(a+b))N+ab\geq N>\lfloor j^2/4\rfloor,
  $$
  so our assumption shows that~(\ref{ineq-1b}) does not hold.

  Assume now that $ab\geq N$. We must then have $j< a+b$ since otherwise
  $$
  (j-(a+b))N+ab\geq ab\geq N>\lfloor j^2/4\rfloor,
  $$
  so that~(\ref{ineq-1b}) does not hold.  We write $a+b=j+k$ with
  $k\geq 1$.

  The left-hand inequality in~(\ref{ineq-1b}) is equivalent to
  $$
  kN\leq ab,
  $$
  and since $ab\leq N^2$, it follows that~$k\leq N$.  Moreover, from
  $a\geq kN/b$, we get
  $$
  a+b\geq \frac{kN}{b}+b\geq 2\sqrt{kN},
  $$
  and therefore
  $$
  k=a+b-j\geq 2\sqrt{kN}-2\sqrt{N+1}.
  $$
  since $j\leq 2\sqrt{N+1}$.

  If $k\geq 5$ and $N>4$, this leads to the contradiction $k>N$ (since
  then $2\sqrt{kN}-2\sqrt{N+1}>\sqrt{kN}$). If $2\leq k\leq 4$, the
  inequality $ 2\sqrt{kN}-2\sqrt{N+1}\leq k$ only happens if $N$ is
  bounded, indeed if $N\leq 8$ by inspection (achieved for $k=2$).

  There only remains the possibility that $k=1$. We then have
  $$
  2\sqrt{N}\leq a+b=j+1\leq 2\sqrt{N+1}+1\leq 2\sqrt{N}+2
  $$
  (for $N\geq 1$).  We write
  $$
  a=\sqrt{N}-h,\quad\quad b=\sqrt{N}+h',
  $$
  where $h\in \Rr$ and $h'\geq 0$.
  The previous result implies
  $$
  2\sqrt{N}\leq 2\sqrt{N}+h'-h\leq 2\sqrt{N}+2
  $$
  hence $0\leq h'-h\leq 2$. But also
  $$
  N>\lfloor j^2/4\rfloor \geq \frac{(a+b-1)^2}{4}-1=
  N+(h'-h-1)\sqrt{N}+\frac{(h'-h-1)^2}{4}-1,
  $$
  which implies that $h'-h-1\leq 1/2$ for $N\geq 8$.\footnote{\ This
    can be improved if we assume that~$N$ is a bit larger.}

  We now get
  $$
  N\leq ab=N+(h'-h)\sqrt{N}-hh'\leq N+\frac{3}{2}\sqrt{N},
  $$
  and
  $$
  (N-a)(N-b)=N^2-(a+b)N+ab=N(N-j)+(ab-N),
  $$
  hence
  $$
  N(N-j)\leq (N-a)(N-b)\leq N(N-j)+\frac{3}{2}\sqrt{N},
  $$
  which gives the desired statement.
\end{proof}

We next look at the interval between $a_j$ and $b_{j-1}$. Given again
arbitrary $(a,b)$ and $j\geq 0$ with $\lfloor j^2/4\rfloor\leq N$, we
have
\begin{equation}\label{ineq-2}
a_j\leq (N-a)(N-b)\leq b_{j-1}
\end{equation}
if and only if 
\begin{equation}\label{ineq-2b}
\lfloor j^2/4\rfloor\leq (j-(a+b))N+ab\leq N.
\end{equation}

\begin{lemma}\label{lm-large}
  Assume that $0\leq a\leq b\leq N$ and that $\lfloor j^2/4\rfloor<N$
  and $\lceil j/2\rceil< N$ so that $a_j\geq
  1$. Then~\emph{(\ref{ineq-2})} holds if and only if
  $$
  (a,b)=(0,j-1)\text{ or } (a,b)=(\lfloor j/2\rfloor,\lceil j/2\rceil).
  $$
\end{lemma}

\begin{proof}
  Suppose first that $ab<N$. The inequality $(N-a)(N-b)\leq b_{j-1}$ is
  equivalent to 
  to $a+b>j-1$ or $a+b=j-1$ and $ab\leq 0$.\footnote{\ Which is
    equivalent to $(a,b)\leqa (0,j-1)$, with notation as in
    Section~\ref{sec-large}.} The second case means that $a=0$, in which
  case $(a,b)=(0,j-1)$. Otherwise, we get $a+b\geq j$. On the other
  hand, if $a+b\geq j+1$, then $(j-(a+b))N+ab\leq -N+an<0$, which
  contradicts~(\ref{ineq-2b}), so only $a+b=j$ remains possible. The
  inequality
  $$
  \lfloor j/2\rfloor\lceil j/2\rceil\leq ab
  $$
  only occurs when $(a,b)=(\lfloor j/2\rfloor,\lceil j/2\rceil)$. So we
  obtain the two allowed cases.

  We suppose now that $ab\geq N$. This implies also that
  $b\geq \sqrt{N}$. The inequality~(\ref{ineq-2b}) implies then that
  $a+b>j$ (the case~$a+b=j$ is impossible, since then
  $ab\leq \lfloor j^2/4\rfloor<N$). We write $a+b=j+k$ with $k\geq
  1$. The left-hand side of~(\ref{ineq-2b}) becomes
  $$
  \lfloor (a+b-k)^2/4\rfloor\leq -kN+ab,
  $$
  hence implies
  $$
  \frac{(a+b)^2}{4}-\frac{k(a+b)}{2}+\frac{k^2}{4}\leq -kN+ab+1.
  $$

  Since $(a+b)^2/4-ab=(a-b)^2/4$, this gives
  $$
  -\frac{k(a+b)}{2}+\frac{k^2}{4}\leq
  \frac{(a-b)^2}{4}-\frac{k(a+b)}{2}+\frac{k^2}{4}\leq -kN+1
  $$
  and this implies that
  $$
  b\geq \frac{a+b}{2}\geq N+\frac{k}{4}-\frac{1}{k}.
  $$

  If $k\geq 3$, this is a contradiction since $b\leq N$; if
  $1\leq k\leq 2$, then it implies that $b\geq N-3/4$, so that $b=N$,
  which is incompatible with~(\ref{ineq-2}) since $a_j\geq 1$.
\end{proof}

We can now show that all the integers in Theorem~\ref{th-table}
occur as primitive gaps in the multiplication table. 

\begin{proposition}\label{pr-existence-gaps}
  Let $N\geq 2$ and let $R$ be the largest non-negative integer such
  that $\lfloor R^2/4\rfloor<N$. The following properties hold.

  \begin{enumth}
  \item For all integers $j$ such that $1\leq j\leq R$, the integer
    $N-\lfloor j^2/4\rfloor$ is a primitive gap in the $N\times N$
    multiplication table, occuring from the interval
    \[
      a_j<b_{j-1}.
    \]
  \item For all integers $j$ such that $1\leq j\leq R-2$, the integer
    $j$ is a primitive gap in the $N\times N$ multiplication table.
  \item If~$R$ is of the form $k^2$, $k^2-1$, $k^2+k$, $k^2+k-1$ for
    some integer~$k\geq 1$, then $R-1$ is a primitive gap in the
    $N\times N$ multiplication table, occuring from the intervals
    \begin{equation}\label{eq-int1}
      (k(k-1))^2=k^2(k^2-2k+1)<(k^2-1)(k^2-2k+2)
    \end{equation}
    for $N=k^2-1$ or $N=k^2$ and
    \begin{equation}\label{eq-int2}
      (k^2+k)(k^2-k)=k^2(k^2-1)<(k^2+k-1)(k^2-k+1)
    \end{equation}
    for $N=k^2+k-1$ or $N=k^2+k$.
  \item If~$R$ is of the form $k^2$ or $k^2+k$ for some
    integer~$k\geq 1$, then~$R$ is a primitive gap in the $N\times N$
    multiplication table occuring from the intervals
    \begin{equation}\label{eq-int3}
      (k^2-k+1)^2<k^2(k^2-2k+3),
    \end{equation}
    if $N=k^2$ and
    \begin{equation}\label{eq-int4}
      k^2(k^2+1)<(k^2+k)(k^2-k+2)
    \end{equation}
    if $N=k^2+k$.
  \end{enumth}
\end{proposition}

\begin{proof}
  We can easily check numerically the statement for small values of~$N$,
  and so we will assume that (say) $N\geq 100$.
  
  From Lemma~\ref{lm-large}, we see that the $N\times N$
  multiplication table contains the primitive gaps
  $N-\lfloor j^2/4 \rfloor$ for $1\leq j\leq R$.
  
  From Lemma~\ref{lm-short} and the fact that
  $$
  N(N-j)\leq (N-1)(N-(j-1))\leq \cdots\leq (N-\lfloor j/2\rfloor)(N-\lceil
  j/2\rceil)
  $$
  with successive gaps
  $$
  (N-(a+1))(N-(j-a-1))-(N-a)(N-(j-a))=j-2a-1
  $$
  for $0\leq a\leq \lfloor j/2\rfloor-1$, we see that the $N\times N$
  multiplication table contains as primitive gaps the integers
  $$
  j-2a-1
  $$
  where, for $1\leq j\leq R$, we let $a$ run over integers such that
  $0\leq a\leq \lfloor j/2\rfloor-1$, \emph{provided}
  $$
  a(j-a)>\frac{3}{2}\sqrt{N}.
  $$
  
  Indeed, this last condition implies that
  $$
  (N-a)(N-(j-a))=N^2-jN+a(j-a)>N(N-j)+\frac{3}{2}\sqrt{N},
  $$
  so that all the products given by the second part of
  Lemma~\ref{lm-short} are at most $(N-a)(N-(j-a))$, implying that the
  gap $j-2a-1$ in the interval
  $$
  (N-a)(N-(j-a))< (N-(a+1))(N-(j-a-1))
  $$
  (which lies between $b_j$ and $a_j$) is primitive.
  
  Note that the values of~$a$ (in the indicated range) with
  $a(j-a)>\frac{3}{2}\sqrt{N}$ is an interval (if the conditions holds
  for~$a$, it also does for~$a+1$, if $a+1$ is still in the allowed
  range).

  Taking $j=R$, we get all gaps $R-2a-1$ for
  $1\leq a\leq \lfloor R/2\rfloor-1$, since $R-1>3\sqrt{N}/2+1$ for
  $N>37$; taking $j=R-1$, we get all gaps $R-2a-2$ for
  $1\leq a\leq \lfloor (R-1)/2\rfloor-1$, since $R-2>3\sqrt{N}/2+1$ for
  $N>64$.  All together, this gives all gaps between $1$ and~$R-3$.
  (If $R$ is odd, one checks that the ``boundaries'' are $2$ and ~$1$,
  respectively, whereas they are $1$ and~$2$ if $R$ is even.)

  

  It remains to show that~$R-2$ is also always a primitive gap, and to
  handle the cases asserting the existence of $R-1$ and $R$ as
  primitive gaps.

  \textbf{Case 1.} We claim that as soon as~$N\geq 8$, the inequality
  $$
  b_{R-1}=N(N-(R-1))<(N-1)(N-(R-2))
  $$
  gives a primitive gap of size~$R-2$ in the multiplication
  table. Indeed, an inequality
  $$
  N(N-R+1)< (N-a)(N-b)< (N-1)(N-(R-2))
  $$
  with~$0\leq a\leq b$ is equivalent to
  $$
  0<N(R-(a+b+1))+ab<R-2.
  $$

  Proceeding as in the proof of Lemma~\ref{lm-short}, we find that
  this implies that $a+b=(R-1)+1=R$ and $N\leq ab$. 
  However, for $a+b=R$, we have $ab\leq \lfloor
  R^2/4\rfloor<N$, so this is impossible.
  
  \textbf{Case 2.} We claim that if $N$ is of one of the forms $k^2$,
  $k^2-1$, $k^2+k$ or~$k^2+k-1$ for $k\geq 4$, then the integer $R-1$
  is a primitive gap given by the interval~(\ref{eq-int1})
  or~(\ref{eq-int2}).  Note that the value of~$R$ is then,
  respectively, equal to $2k-1$, $2k-1$, $2k$, $2k$. It suffices then
  to show that the following intervals do indeed define primitive
  gaps:
  \begin{align*}
    (k(k-1))^2=k^2(k^2-2k+1)<(k^2-1)(k^2-2k+2),&\quad N=k^2\\
    (k^2+k)(k^2-k)=k^2(k^2-1)<(k^2+k-1)(k^2-k+1),&\quad N=k^2+k
  \end{align*}

  Indeed, the first one takes care of the cases $N=k^2$ or $k^2-1$,
  since the elements on each side of the inequality belong to the
  $(k^2-1)$-table, and the second one of the cases $N=k^2+k$ or
  $k^2+k-1$.
  
  We can handle both inequalities together, since they are the same as
  $$
  N(N-R)<(N-1)(N-R+1).
  $$

  As before, an inequality
  $$
  N(N-R)< (N-a)(N-b)< (N-1)(N-R+1)
  $$
  with~$0\leq a\leq b$ is equivalent to
  $$
  0<N(R-(a+b)N)+ab<R-1,
  $$
  and is only possible if $a+b=R+1$. But then
  $ab\leq \lfloor (R+1)^2/4\rfloor$, and by inspection we have
  $\lfloor (R+1)^2/4\rfloor=N$ in the present case (we have
  $(R+1)^2/4=k^2$ if $N=k^2$ and $(R+1)^2/4=k^2+k+1/4$ if $N=k^2+k$).

  \textbf{Case 3.} We claim that if $N$ is of one of the forms $k^2$ or
  $k^2+k$,
  then the integer $R$ is a primitive gap given by the
  interval~(\ref{eq-int3}) or~(\ref{eq-int4}). Note that the value
  of~$R$ is then, respectively, equal to $2k-1$ and $2k$. It suffices
  then to show that the following intervals do indeed define primitive
  gaps:
  \begin{align*}
    (k^2-k+1)^2<k^2(k^2-2k+3),&\quad N=k^2,\\
    k^2(k^2+1)<(k^2+k)(k^2-k+2),&\quad N=k^2+k.
  \end{align*}

  If we transcribe in terms of~$N$ and~$R$ these inequalities, and
  consider when a product $(N-a)(N-b)$ lies strictly between these
  quantities, we find that this is equivalent to
  \begin{gather*}
    -R<N(R-(a+b+2))+ab<0\\
    -R<N(R-(a+b+2))+ab<0,
  \end{gather*}
  respectively. (Each of these depends on the precise relation between
  $N$ and~$R$; for instance, for the first of these, we use the fact
  that
  $$
  N+\Bigl(\frac{R-1}{2}\Bigr)^2=2N-R
  $$
  when~$N=k^2$ and~$R=2k-1$.) In particular, there is really only one
  case.
  Arguing as previously,\footnote{\ It is helpful to note here that
    $\sqrt{kN-R}+2\geq \sqrt{kN}$ as soon as $N\geq 10$, for instance.}
  we find that $R-(a+b+2)=-1$, hence $a+b=R-1$, and the inequalities are
  equivalent to $N-R<ab<N$. But for $a+b=R-1$, we get
  $$
  ab\leq \Bigl\lfloor\frac{(R-1)^2}{4}\Bigr\rfloor,
  $$
  which in both cases here is equal to $N-R$. (For instance, if
  $N=k^2$, $R=2k-1$, we get $(R-1)^2/4=(k-1)^2=N-R$.)
\end{proof}

We record another corollary of the results of this section which will be
useful in the next section.

\begin{corollary}\label{cor-large-range}
  Let $N\geq 2$ and let $R$ be the largest non-negative integer such
  that $\lfloor R^2/4\rfloor<N$. All the primitive gaps in the
  multiplication  table in the interval
  \[
    b_{R-2}=\bj{(R-2)}<N^2
  \]
  are either $\leq R-3$ or in $\mathcal{B}_N$. In the second case, these
  gaps are isolated.
\end{corollary}

\begin{proof}
  The indicated range is split in the subintervals
  \[
    b_{R-2}<a_{R-2}<b_{R-3}<\cdots < b_2<a_2<b_1<N^2.
  \]

  The gaps $b_{j-1}-a_j$ with $1\leq j\leq R-2$ are primitive and belong
  to $\mathcal{B}_N$ by Proposition~\ref{pr-existence-gaps},
  (1). Moreover, for $2\leq j\leq R-2$, we have the further splitting
  \[
    b_j=N(N-j)\leq (N-1)(N-(j-1))\leq \cdots\leq \aj{j}=a_j,
  \]
  with successive gaps $j-2a-1$ for $0\leq a\leq \lfloor j/2\rfloor-1$,
  so that the primitive gaps within the interval $b_j<a_j$ are all
  $\leq j-1\leq R-3$.

  Finally, again from the previous results, any gap in $\mathcal{B}_N$
  occurs uniquely in the indicated range.
\end{proof}

\section{The list of gaps is exhaustive}

In this section, we conclude the proof of Theorem~\ref{th-table} by
showing that there are no other gaps than those described in the
statement.

We begin with a lemma which will be used to identify the special forms
of the values of~$N$ for which either~$R-1$ or~$R$ is a primitive gap
in~$\mathcal{S}_N$.

\begin{lemma}\label{lm-r}
  Let~$k\geq 1$ be an integer.
  \begin{enumth}
  \item Let~$N=k^2+j$ for some integer~$j$ with $1\leq j\leq k$. We
    then have $R=2k$.
  \item Let~$N=k^2+k+j$ for some integer~$j$ with $1\leq j\leq k+1$. We
    then have $R=2k+1$.
  \item We have
    \[
      \Bigl\lfloor \frac{(R+1)^2}{4} \Bigr\rfloor =
      \begin{cases}
        N+1 & \text{ if  $N$ is of the form $k^2-1$ or $k^2+k-1$},
        \\
        N+2 &\text{ if $N$ is of the form $k^2-2$ or $k^2+k-2$},
        \\
        \geq N+3&\text{ otherwise.}
      \end{cases}
    \]
  \end{enumth}
\end{lemma}

\begin{proof}
  The first two statements are elementary.
  Suppose then that $N=k^2+j$ with $1\leq j\leq k$. Then $R=2k$, hence
  \[
    \Bigl\lfloor \frac{(R+1)^2}{4} \Bigr\rfloor =k^2+k,
  \]
  which is equal to $N+1$ if $j=k-1$, i.e. $N=k^2+k-1$, to $N+2$ if
  $N=k^2+k-2$ and otherwise is $\geq N+3$. The case where $N=k^2+k+j$
  with $1\leq j\leq k$, and $R=2k+1$, is handled similarly.
\end{proof}

We now start the proof of Theorem~\ref{th-table}. We will proceed by
induction on~$N$. The key idea is that most elements of the
$(N+1)\times (N+1)$ multiplication table are present in the $N \times N$
table, so in the range where these gaps are quite structured, we only
need keep track of how the new entries interact with the old. As is
often the case, the argument is easier if we use a stronger induction
hypothesis than the statement of Theorem~\ref{th-table} for a
given~$N$. Precisely, the assumption (denoted $\textup{H}_N$) will be
\begin{center}
  $\textup{H}_N$: \textit{The set of gaps coincides with the union
    $\mathcal{S}_N\cup \mathcal{B}_N$ of Theorem~\ref{th-table}, and
    moreover the gaps of the form $N-\lfloor j^2/4\rfloor$, as well as
    those equal to~$R-1$ or~$R$ (when they exist) are isolated.}
\end{center}

Note that $\textup{H}_N$ implies the statement of
Theorem~\ref{th-table} for a given~$N$. We will prove that
$\textup{H}_N$ holds for $N\geq 2$ by induction on~$N$.  For any small
enough set of values of~$N$, we can check the statement $\textup{H}_N$
by hand or by computer. We do this for $N\leq 6$, using the
corresponding lists of primitive gaps with multiplicities:
\begin{gather*}
  (1,1,1)\text{ for } N=2,\quad (3,3,1,1)\text{ for } N=3,\quad
  (7,4,2,1,1)\text{ for } N=4,
  \\
  (11,8,2,1,1,1)\text{ for } N=5,\quad (18,9,4,1,1,1,1)\text{ for } N=6,
\end{gather*}
with notation as in Example~\ref{ex-examples}.

We assume now that~$\textup{H}_N$ holds for some integer $N\geq 5$,
and will prove that $\textup{H}_{N+1}$ also holds. We
then have $R\geq 4$. We denote by $T$ the maximum of the
set~$\mathcal{S}_N$ and by $R'$ the integer~$R$ for the
$(N+1)\times (N+1)$-multiplication table.  For conciseness and to
avoid any confusion, we speak of $N$-gap (resp. $(N+1)$-gap) for a
primitive gap in the $N\times N$ multiplication table (resp
$(N+1)\times (N+1)$ multiplication table).  For simplicity, we will
also say that \emph{there is an exceptional $N$-gap} if~$T>R-2$.

The induction step will be split in three subcases:
\begin{enumerate}
\item[] Case 1: there is no exceptional $N$-gap and $R'=R$,
\item[] Case 2: there is an exceptional $N$-gap and $R'=R$,
\item[] Case 3: there is an exceptional $N$-gap and $R'=R+1$.
\end{enumerate}

This is indeed an exhaustive list of possibilities, since
$R\leq R'\leq R+1$ and since $R'=R+1$ only occurs when~$N$ is of the
form $k^2$ or $k^2+1$ (from Lemma~\ref{lm-r}), in which case there is
an exceptional $N$-gap.

\par
\textbf{Case 1:} There is no exceptional $N$-gap and $R'=R$.
\par

By Corollary~\ref{cor-large-range}, applied to $N+1$, we know that all
$(N+1)$-gaps in the range
\[
  (N+1)(N+1-(R'-2))=(N+1)(N-R+3)<(N+1)^2
\]
are either $\leq R'-2$ or are gaps in~$\mathcal{B}_{N+1}$ which occur
only once.  In addition, the induction assumption $\textup{H}_N$
implies that all $N$-gaps located in the interval
\[
  1<N(N-R+2)
\]
are $\leq R-2$ (indeed, $\textup{H}_N$ implies that any larger $N$-gap
is isolated, and inspection of their unique appearance, given by
Proposition~\ref{pr-existence-gaps}, (3) and~(4) and
Corollary~\ref{cor-large-range}, imply this property). Any $(N+1)$-gap
which lies in this interval is therefore also $\leq R-2=R'-2$.

We now claim that the interval
\[
  N(N-R+2)<(N+1)(N-R+3)
\]
contains no $(N+1)$-gap $\geq R-1$, except for an isolated exceptional
$(N+1)$-gap equal to~$R'-1$ if~$N+1$ is of the form $k^2-1$ or
$k^2+k-1$. This will conclude the proof of $\textup{H}_{N+1}$.

The upper bound in the range above satisfies the inequality
\[
  (N+1)(N-R+3)<N^2
\]
since $R\geq 4$ (which follows from our assumption that $N\geq 5$).
The lower bound is equal to $b_{R-2}$, so the range in question is
within the range where Corollary~\ref{cor-large-range} applies to the
$N\times N$ table. Using again $R=R'$, this corollary implies that any
$(N+1)$-gap $\geq R'-1$ must lie within an interval $a_j<b_{j-1}$ of
the $N\times N$-table. Since
\[
  b_{R-4}=N(N-R+4)>(N+1)(N-R+3)
\]
(again because $R\geq 4$), the only possible intervals are
\[
  a_{R-2}<b_{R-3}\quad\text{and}\quad a_{R-3}<b_{R-4}.
\]

For the former, observe that
\[
  a_{R-2}< (N+1)(N-R+2)< b_{R-3},
\]
since the differences are
\begin{gather*}
  b_{R-3}-(N+1)(N-R+2)=R-2
  \\
  (N+1)(N-R+2)-a_{R-2}= N-R+2-\Bigl\lfloor \frac{(R-2)^2}{4}
  \Bigr\rfloor=N+1-\Bigl\lfloor
  \frac{R^2}{4}\Bigr\rfloor\geq 1.
\end{gather*}

In fact, the last difference satisfies
\[
  N+1-\Bigl\lfloor \frac{R^2}{4}\Bigr\rfloor\leq R-2
\]
as soon as $N\geq 7$ (this can be checked using
Lemma~\ref{lm-r}\footnote{\ If $N=k^2+j$ with $1\leq j\leq k$, for
  instance, then the difference is $j+1\leq k+1$ and $R-2=2k-2$, so
  $k\geq 3$ suffices.}).
Thus, the gaps in the $(N+1)\times (N+1)$ table between $a_{R-2}$ and
$b_{R-3}$ are $\leq R-2$.

For the second interval $a_{R-3}<b_{R-4}$, we observe that
\[
  a_{R-3}<(N+1)(N-R+3),
\]
since the difference is
\[
  N+R+1-\Bigl\lfloor\frac{(R+1)^2}{4}\Bigr\rfloor= R-2+
  \Bigl(N+3-\Bigl\lfloor\frac{(R+1)^2}{4}\Bigr\rfloor\Bigr),
\]
and the second term in the right-hand side is $\leq 2$ by
Lemma~\ref{lm-r}, (3).  Thus the largest possible $(N+1)$-gap in the
allowed range between $a_{R-3}$ and $b_{R-4}$ is the quantity
$R-2+(N+3-S)$ where
\[
  S=\Bigl\lfloor \frac{(R+1)^2}{4}\Bigr\rfloor.
\]

If~$S\geq N+3$, then any occuring gap is $\leq R-2$. On the other
hand, by Lemma~\ref{lm-r}, (3), we have $S=N+1$ if and only if $N$ is
of the form $k^2-1$ or $k^2+k-1$. But it this were the case, there
would be an exceptional $N$-gap, contrary to our assumption in this
subcase of the induction. Thus the remaining possibility is $S=N+2$,
which by loc. cit. occurs if and only if $N$ is of the form $k^2-2$ or
$k^2+k-2$.

In this case, the interval
\[
  a_{R-3}<(N+1)(N-R+3)
\]
coincides with
\[
  (k(k-1))^2<(k^2-1)(k^2-2k+2).
\]

Proposition~\ref{pr-existence-gaps} shows that this interval does define
a primitive exceptional gap of size $R-1=R'-1$ in the
$(N+1)\times (N+1)$-table.  This concludes the proof of the claim, and
hence of Case~1.


\par
\textbf{Case 2:} There is an exceptional $N$-gap and $R'=R$.
\par
These conditions imply that $N$ is of the form $k^2-1$ (with an
exceptional $N$-gap equal to $R-1$) or $k^2+k-1$ (with exceptional
$N$-gaps equal to $R-1$ or~$R$).

To prove that all $(N+1)$-gaps are in $\mathcal{S}_{N+1}$
or~$\mathcal{B}_{N+1}$, we can argue as the beginning of Case~1
(since the condition $R'=R$ was also satisfied there). We find, using
induction, that the only non-exceptional $(N+1)$-gaps of size
$\geq R'-2$ which could exist outside of $\mathcal{B}_{N+1}$ would be
in the interval
\[
  a_{R-3}<(N+1)(N-R+3),
\]
or would be given by an exceptional $N$-gap of size $R-1$ from
$\mathcal{E}_N$.

Using the values of $N$ and~$R$ above, we find that satisfies
\[
  (N+1)(N-R+3)-a_{R-3}=2k-1=R'
\]
(use Lemma~\ref{lm-r}) and that, for $N=k^2-1$ as well as $N=k^2+k-1$,
this interval coincides exactly with the interval for the primitive gap
of size $R'$ described in Corollary~\ref{cor-large-range}.

Furthermore, by the induction assumption $\mathrm{H}_N$, there is a
unique exceptional $N$-gap of size $R-1$, given by the interval
\[
  (k(k-1))^2=k^2(k^2-2k+1)<(k^2-1)(k^2-2k+2)
\]
for $N=k^2-1$ or
\[
  (k^2+k)(k^2-k)=k^2(k^2-1)<(k^2+k-1)(k^2-k+1)
\]
for $N=k^2+k-1$. By Proposition~\ref{pr-existence-gaps}, (3), this is
also a primitive gap, of the same size, in the $(N+1)\times (N+1)$
multiplication table.  We conclude from these arguments that the first
part of~$\mathrm{H}_{N+1}$ holds.

For the other assertions, it suffices to prove that there is no gap of
size $\geq R'$ in the range
\[
  1\leq N(N-R'+2),
\]
according to Corollary~\ref{cor-large-range}, and this holds again by
induction, since $R'=R$.

Finally, the reasoning used at the beginning of Case~2 shows that the
exceptional $(N+1)$-gaps are isolated.  Thus we deduce
that~$\mathcal{H}_{N+1}$ holds.

\par
\textbf{Case 3:} There is an exceptional $N$-gap and $R'=R+1$.
\par
These conditions imply that $N$ is of the form $k^2$ (with $R=2k-1$,
$R'=2k$) or $k^2+k$ (with $R=2k$, $R'=2k+1$). Moreover, both $R-1$
and~$R$ are exceptional~$N$-gaps, and there is no exceptional
$(N+1)$-gap.

Arguing as in Case~1, we see that any primitive $(N+1)$-gap of size
$\geq R'-2=R-1$ is contained in the interval
\[
  N(N-R+2)<(N+1)(N-R'+3)=(N+1)(N-R+2)
\]
or within the interval of an exceptional $N$-gap.

We consider the first possibility. The interval is again within the
range of Corollary~\ref{cor-large-range} for the
$N\times N$-multiplication table, so the only possible primitive gaps
of size $\geq R'-2$ must be contained in one of the gaps $b_{j-1}-a_j$
in $\mathcal{B}_N$. However, we see that the only possibility~is
\[
  a_{R-2}<b_{R-3}
\]
because $b_{R-3}>(N+1)(N-R+2)$ if $R\geq 3$. However, we compute in
both cases that
\[
  (N+1)(N-R+2)-a_{R-2}=k+1\leq R-2.
\]

We now consider the exceptional $N$-gaps. Since $R-1=R'-2$, the gap of
size $R-1$ is not an issue. On the other hand, the exceptional $N$-gap
of size $R=R'-1$ is isolated (using $\mathrm{H}_N$) and corresponds to
the interval
\begin{gather*}
  (k^2-k+1)^2<k^2(k^2-2k+3)\quad\quad\text{ if } N=k^2\\
  k^2(k^2+1)<(k^2+k)(k^2-k+2)\quad\quad\text{ if } N=k^2+k,
\end{gather*}
by the induction hypothesis. If $N=k^2$, we observe that
\[
  (k^2-k+1)^2<(k+1)^2(k^2-2k+2)<k^2(k^2-2k+3)
\]
with differences $2k-2=R'-2$ and $1$; if $N=k^2+k$, we observe that
\[
  k^2(k^2+1)<
  (k^2+k+1)(k^2-k+1)
  <(k^2+k)(k^2-k+2)
\]
with differences $2k-1=R'-2$ and $1$. In both cases, there is
therefore no gap of size $R'-1$ in the above intervals, and we
conclude that the first part of $\mathrm{H}_{N+1}$ holds.

The above has also shown that the gaps in $\mathcal{B}_{N+1}$ are
those from Proposition~\ref{pr-existence-gaps}, and hence are
isolated.  Since there are no exceptional $(N+1)$-gap, this shows that
$\mathrm{H}_{N+1}$ holds.

This concludes the proof of Case~3, and hence finally of
Theorem~\ref{th-table}.


\begin{thebibliography}{CCC}

\bibitem{b-k} R.P. Brent and H.T. Kung: \textit{The area-time
    complexity of binary multiplication}, J. Assoc. Comput. Mach.
  28 (1981), 521--534.

\bibitem{bppw} R. Brent, C. Pomerance, D. Purdum and J. Webster:
  \textit{Algorithms for the multiplication table problem}, Integers
  21 (2021), Paper No. A92, 19.
  
\bibitem{efg} S. Eberhard, K. Ford and B. Green:
  \textit{Permutations fixing a $k$-set}, IMRN (2016), 6713--6731.

\bibitem{edwards} H.M. Edwards: \textit{Riemann's zeta function},
  Academic Press, 1974.
  
\bibitem{ford} K. Ford: \textit{The distribution of integers with a
    divisor in a given interval}, Annals of Math. 168 (2008), 367--433.
  
\bibitem{ganguly-guria} S. Ganguly and R. Guria: \textit{Lattice points
    on determinant surfaces and the spectrum of the automorphic
    laplacian}, preprint (2024), \url{arXiv:241004637}. 

\bibitem{gk} S.W. Graham and G. Kolesnik: \textit{van der Corput's
    method of exponential sums}, L.M.S Lecture Note Series 126,
  Cambridge Univ. Press, 1991.

\bibitem{granville} A. Granville: \textit{Harald Cramér and the
    distribution of prime numbers}, Scand. Actuarial J. (1995),
  12--28.
  
\bibitem{granville2} A. Granville: \textit{The anatomy of integers and
    permutations}, preprint,
  \url{https://www.dms.umontreal.ca/~andrew/PDF/AnatomyForTheBook.pdf}


\bibitem{katz-sarnak} N.M. Katz and P. Sarnak: \textit{Zeros of zeta
    functions and symmetries}, Bull. of the AMS 36 (1999), 1--26.
  
\bibitem{k-t} D. Koukoulopoulos and J. Thiel: \textit{Arrangements of
    stars on the American flag}, Amer. Math. Monthly 119 (2012),
  443--450.
 
\bibitem{pnt} E. Kowalski: \textit{An introduction to probabilistic
    number theory}, Cambridge Studies in Advanced Mathematics 192,
  Cambridge University Press, 2021.

\bibitem{m-s} J. Marklof and A. Strömbergsson: \textit{The three gap
    theorem and the space of lattices}, American Math. Monthly 124
  (2017), 741--745.

\bibitem{maynard} J. Maynard: \textit{Small gaps between primes},
  Annals of Math. 181 (2015), 383--413.

\bibitem{montgomery} H.L. Montgomery: \textit{The pair-correlation of
    zeros of the zeta function}, in Analytic Number Theory,
  Proc. Sympos. Pure Math. XXIV, AMS, 1972, 181--193.
  
\bibitem{parigp} The PARI~Group, PARI/GP version \texttt{2.14.0},
  Université de Bordeaux, 2022, \url{https://pari.math.u-bordeaux.fr/}.

\bibitem{sarnak} P. Sarnak: \textit{Spectra and eigenfunctions of
    Laplacians}, in Partial differential equations and their
  applications, CRM Proc. Lecture Notes 12, AMS, 1997, 261--276.
  
\bibitem{sos} V. Sós: \textit{On the theory of diophantine
    approximatelys I}, Acta Math. Acad. Sci. Hungar. 8 (1957),
  461--472.

\bibitem{tenenbaum} G. Tenenbaum: \textit{Sur la probabilité qu'un
    entier possède un diviseur dans un intervalle donné}, Compositio
  Math. 51 (1984), 243--263.


\bibitem{tschinkel}
 Y. Tschinkel: \textit{About the cover: on the distribution of
   primes--Gauss’ tables}, Bulletin of the AMS 43 (2006), 89--91.

\bibitem{zagier}
 D. Zagier: \textit{A one-sentence proof that every
    prime $p\equiv 1\bmod{4}$ is a sum of two squares}, American
  Math. Monthly 97 (1990), 144.

\bibitem{zhang}
 Y. Zhang: \textit{Bounded gaps between primes}, Annals
 of Math 179 (2014), 1121--1174.
  
\end{thebibliography}
\end{document}